\documentclass{amsart}
\usepackage{tipa}
\usepackage{amssymb,amsmath}
\usepackage{stmaryrd}

\usepackage[cmtip,all]{xy}
\usepackage{graphicx,subfigure}

\newtheorem{theorem}{Theorem}[section]

\theoremstyle{definition}

\newtheorem{corollary}[theorem]{Corollary}
\newtheorem{remark}[theorem]{Remark}
\newtheorem{question}[theorem]{Question}

\numberwithin{equation}{section}

\begin{document}
\title[Integers representable as the product]{Integers representable as the product of the sum of four integers with the sum of their reciprocals}
\author{YONG ZHANG}

\address{School of Mathematics and Statistics, Changsha University of Science and Technology,
Changsha 410114, People's Republic of China}

\address{School of Mathematical Sciences, Zhejiang University, Hangzhou 310027, People's Republic of China}
 \email{zhangyongzju$@$163.com}

\thanks{This research was supported by the National Natural Science Foundation of China (Grant No.~11501052).}

\subjclass[2010]{Primary 11D25; Secondary 11D72, 11G05}

\date{}

\keywords{the sum of integers, the sum of integers' reciprocals,
elliptic curves}

\begin{abstract}
By the theory of elliptic curves, we study the integers
representable as the product of the sum of four integers with the
sum of their reciprocals and give a sufficient condition for the
integers with a positive representation.
\end{abstract}
\maketitle

\section{Introduction}
There are many interesting results about the representation of
integers $n$ by special symmetric polynomials $f(x,y,z)$, such
as\[\begin{split}
n=f(x,y,z)=&(x+y+z)\bigg(\frac{1}{x}+\frac{1}{y}+\frac{1}{z}\bigg),\frac{(x+y+z)^3}{xyz},\\
&\frac{x}{y}+\frac{y}{z}+\frac{z}{x},\frac{x}{y+z}+\frac{y}{z+x}+\frac{z}{x+y},xy+yz+zx.\end{split}\]
We refer to
\cite{Bremner-Guy-Nowakowski,Brueggeman,Bremner-Guy,Bremner-Macleod,Cai}.

For $m\geq2$, it is interesting to find nonzero integers $x_1,\cdot
\cdot \cdot,x_m$ such that
\begin{equation}
n=(x_1+\cdot \cdot \cdot+x_m)\bigg(\frac{1}{x_1}+\cdot \cdot \cdot+\frac{1}{x_m}\bigg).
\end{equation}
If $x_i,i=1,..,m$ are all positive, we say $n$ has a positive
representation.

When $m=2$, Eq. (1.1) has integer solutions only for $n=0,4$.

When $m=3$, Eq. (1.1) reduces to
\begin{equation}
n=(x+y+z)\bigg(\frac{1}{x}+\frac{1}{y}+\frac{1}{z}\bigg),
\end{equation}
which is known as Melvyn Knight's problem: which integers $n$ can be
represented as the product of the sum of three integers with the sum
of their reciprocals. In 1993, A. Bremner, R.K. Guy and R.J.
Nowakowski \cite{Bremner-Guy-Nowakowski} studied this problem and
gave the necessary and sufficient conditions for the integers with a
positive representation. In 2010, R. Kozuma \cite{Kozuma} got two
sufficient conditions for an integer $n$ which cannot be expressed
in the form of Eq. (1.2) by considering the elliptic curve $E_n:~y^2
+(n-3)xy+(n-1)y=x^3$ over $\mathbb{Q}$.

When $m=4$, the authors of \cite{Bremner-Guy-Nowakowski} gave two
remarks about Eq. (1.1): (1) ``If
\begin{equation}
n=(x+y+z+w)\bigg(\frac{1}{x}+\frac{1}{y}+\frac{1}{z}+\frac{1}{w}\bigg),
\end{equation}
then there are infinitely many representations of every integer $n$,
e.g.,
\[(x,y,z,w)=(m^2+m+1,m(m+1)(n-1),(m+1)(n-1),-m(n-1))\]
for any integer $m$"; (2) ``The least $n$ with a positive
representation is 16, and it seems likely that for $n\geq16$ there
is always such a representation". To prove the second one is a
difficult problem.

In this paper, we mainly consider the positive integer solutions of
Eq. (1.3), i.e., the positive representation of $n$. The right side
of Eq. (1.3) can be written as
\[16+\frac{(x-y)^2}{xy}+\frac{(x-z)^2}{xz}+\frac{(x-w)^2}{xw}+\frac{(y-z)^2}{yz}+\frac{(y-w)^2}{yw}+\frac{(z-w)^2}{zw},\]
so if $x,y,z,w$ are all positive integers, no integer less than 16
can be represented, and 16 only with $x=y=z=w.$ In view of the
homogeneity of Eq. (1.3), it only needs to study the positive
rational solutions $(x,y,z)$ of the equation
\begin{equation}
n=(x+y+z+1)\bigg(\frac{1}{x}+\frac{1}{y}+\frac{1}{z}+1\bigg)
\end{equation}
for $n>16.$ By the theory of elliptic curves, we give a sufficient
condition for the integers $n$ with a positive representation.

\begin{theorem} For positive integers $n>16$ and positive rational numbers $z$ with $nz-(z+1)^2>0$, if the elliptic curve
\[E_{n,z}:~Y^2=X(X^2+(nz(nz-2z^2-8z-2)+(z^2-1)^2)X+16nz^3(z+1)^2)\]
has a rational point $P=(X,Y)$ satisfying the condition
\begin{equation}
X<0,-\frac{X(X-2z(nz+(z+1)^2))}{2z}<Y<((z+1)^2-nz)X,
\end{equation}
then Eq. (1.4) has a positive rational solution $(x,y,z)$, i.e., $n$
can be represented as the product of the sum of four
\textbf{positive} integers with the sum of their reciprocals.
\end{theorem}

\begin{remark}
(1) If Eq. (1.3) has a positive integer solution $(x,y,z,z)$, then
\begin{equation}
n=(x+y+2z)\bigg(\frac{1}{x}+\frac{1}{y}+\frac{2}{z}\bigg),
\end{equation}
and $E_{n,z}$ reduces to
\[E_{n,1}:~Y^2=X(X^2+(n^2-12n)X+64n).\]
We can discuss $E_{n,1}$ by the same methods of
\cite{Bremner-Guy-Nowakowski}. For integer $k$, Eq. (1.6) has
positive integer solutions:
\[
n=4L_{4k}+17,(x,y,z)=(F_{2k-1},F_{2k+1},2F_{2k-1}L_{2k}F_{2k+1}),k\geq
1,
\]
where $F_{k},L_{k}$ are the Fibonacci and Lucas numbers, satisfying
$u_{k+1}=u_k+u_{k-1},F_1=F_2=L_1=1,L_2=3.$

(2) If Eq. (1.3) has a positive integer solution $(x,x,y,y)$, then
\[n=4(x+y)\bigg(\frac{1}{x}+\frac{1}{y}\bigg).\]
Solve it for $x$, we have
\[x=\bigg(\frac{n-8}{8}\pm\frac{\sqrt{n^2-16n}}{8}\bigg)y.\]
In view of $n^2-16n$ is a square and $n>16,x,y>0$, then $n=18,25.$
So 18 and 25 are the only positive integers with a positive
representation $(x,x,y,y)$.

(3) If Eq. (1.3) has a positive integer solution $(x,y,y,y)$, by the
same method of (2), 20 is the only positive integer with a positive
representation $(x,y,y,y)$.
\end{remark}

From (1) of Remark 1.2, we get
\[\begin{split}
4L_{4k}+17=&(F_{2k-1}+F_{2k+1}+2F_{2k-1}L_{2k}F_{2k+1}+2F_{2k-1}L_{2k}F_{2k+1})\\
&\times
\bigg(\frac{1}{F_{2k-1}}+\frac{1}{F_{2k+1}}+\frac{1}{2F_{2k-1}L_{2k}F_{2k+1}}+\frac{1}{2F_{2k-1}L_{2k}F_{2k+1}}\bigg).
\end{split}\]
This leads to the following interesting result.

\begin{corollary} There are infinitely many positive integers $n$ can be representable as the product of the sum of four \textbf{positive} integers with the
sum of their reciprocals.
\end{corollary}

\vskip10pt
\section{The proof of theorem 1.1}
\emph{Proof of Theorem 1.1.} From Eq. (1.4), we have
\[(yz+y+z)x^2+((1+z)y^2+(z^2+4z+1-nz)y+z^2+z)x+yz(y+z+1)=0.\]
Consider the above equation as a quadratic equation of $x$, if it
has rational solutions, then the discriminant $\Delta(x)$ should be
a perfect square. So we need to study the rational points $(y,t)$ on
the quartic elliptic curve
\begin{equation}
\begin{split}
t^2=\Delta(x)=&(z+1)^2y^4+2(1+z)((z+1)^2-nz)y^3\\
&+(z^2n^2-2z(z^2+4z+1)n+(z^2+4z+1)(z+1)^2)y^2\\
&+2z(1+z)((z+1)^2-nz)y+z^2(z+1)^2.
\end{split}
\end{equation}
In view of a method of Mordell (see \cite[p.77]{Mordell} or
\cite[p.476, Proposition 7.2.1]{Cohen1}), $(2.1)$ is birationally
equivalent to the elliptic curve
\[
E_{n,z}:~Y^2=X(X^2+(nz(nz-2z^2-8z-2)+(z^2-1)^2)X+16nz^3(z+1)^2),
\]
by the following mapping
\begin{equation}
\begin{cases}
\begin{split}
y=&\frac{Y+X(nz-(z+1)^2)}{2(z+1)(X-4nz^2)},\\
t=&(8nz^2(nz-(z+1)^2)Y-X^3+12nz^2X^2\\
&+8nz^2(nz(nz-2z^2-8z-2)+(z^2+1)(z+1)^2)X\\
&+64n^2z^5(z+1)^2)/(4(z+1)(X-4nz^2)^2),
\end{split}
\end{cases}
\end{equation}
and its inverse transformation is
\[
\begin{cases}
\begin{split}
X=&-2(z+1)((-z-1)y^2+(nz-(z+1)^2)y+t-z^2-z),\\
Y=&2(z+1)(2(z+1)^2y^3-3(z+1)(nz-(z+1)^2)y^2\\
&+(nz(nz-2z^2-8z-2)+(z+1)(z^3+5z^2+5z+1-2t))y\\
&+(t-z^2-z)(nz-(z+1)^2)).
\end{split}
\end{cases}
\]

The discriminant of $E_{n,z}$ is
\[\Delta(n,z)=256(z+1)^4z^6(z^2n^2-2z(z^2+6z+1)n+(z-1)^4)(nz-(z+1)^2)^2n^2.\]
For $z>0$, $\Delta(n,z)=0$ if and only if $n=0,4,16$. So $E_{n,z}$
is nonsingular for $n>16$ and $z>0$.

It is easy to check that the point
\[P=\big(4z(1+z)^2,4z(1+z)^2(nz-(z+1)^2)\big)\]
lies on $E_{n,z}.$

By the group law, we obtain
\[
\begin{split}
[2]P=&\big(4z^2,-4z^2(nz+z^2+1)\big),\\
[4]P=&\bigg(\frac{4z^2(n(z+1)^2-z)^2}{(nz+z^2+1)^2},\frac{4z^2(n(z+1)^2-z)}{(nz+z^2+1)^3}(z^2(1+z)^2n^3-z^2(4z^2\\
&+7z+4)n^2+(z^6+2z^5+9z^4+12z^3+9z^2+2z+1)n+z^5+z)\bigg).
\end{split}
\]

Let $E_{n,1}$ be the specialization of $E_{n,z}$ at $z=1$, and the
specialization of point $[4]P$ at $z=1$ is
\[[4]P_1=\bigg(\frac{4(4n-1)^2}{(n+2)^2},\frac{4(4n-1)(4n^3-15n^2+36n+2)}{(n+2)^3}\bigg).\]
When the numerator of the $x$-coordinate of $[4]P_1$ is divided by
$(n+2)^2$, the remainder equals $r=-288n-252$ and $r\neq 0$ when
$n>16,$ so the $x$-coordinate of $[4]P_1$ is not a polynomial. For
$17\leq n\leq284$ one can check that $r/(n+2)^2$ is not an integer,
and that it is nonzero and less than 1 in modulus for $n>285$. Hence
for all $n>16$ the point $[4]P_1$ has nonintegral $x$-coordinate and
hence, by the Nagell-Lutz Theorem (see \cite[p.56]{Silverman-Tate}),
is of infinite order. Then $E_{n,z}$ has positive rank in the field
$\mathbb{Q}(z)$. According to the Specialization Theorem of
Silverman (see \cite[p.457, Theorem 20.3]{Silverman}), for all but
finite many positive rational numbers $z$, $E_{n,z}$ has positive
rank.

From Mazur's Theorem (see \cite[p.58]{Silverman-Tate}), a point of
finite order on an elliptic curve has maximal order 12. Let us
compute the expression $[m]P=(X(m),Y(m))$ with $1\leq m\leq 12,$ and
determine those $z>0$ such that denominator of $X(m)$ has a zero at
$z$. By some calculations, we find that for every $z>0$, $[m]P$ is
not a point of finite order, so $P$ is a point of infinite order.
Then $E_{n,z}$ has positive rank over $\mathbb{Q}$ for $z>0,n>16$.

From the transformation $(2.2)$, we get
\[\begin{split}
x=-\frac{X^2-2z(nz+(z+1)^2)X+2zY}{2(1+z)((nz-z^2-1)X-8nz^3+Y)},y=\frac{Y+X(nz-(z+1)^2)}{2(z+1)(X-4nz^2)}.
\end{split}\]
In view of $x>0,y>0$, we have four cases:
\[\begin{split}
(1)~&X^2-2z(nz+(z+1)^2)X+2zY>0,(nz-z^2-1)X-8nz^3+Y<0;\\
&Y+X(nz-(z+1)^2)>0,X-4nz^2>0; \\
or~(2)~&X^2-2z(nz+(z+1)^2)X+2zY>0,(nz-z^2-1)X-8nz^3+Y<0;\\
&Y+X(nz-(z+1)^2)<0,X-4nz^2<0;\\
or~(3)~&X^2-2z(nz+(z+1)^2)X+2zY<0,(nz-z^2-1)X-8nz^3+Y>0;\\
&Y+X(nz-(z+1)^2)>0,X-4nz^2>0;\\
or~(4)~&X^2-2z(nz+(z+1)^2)X+2zY<0,(nz-z^2-1)X-8nz^3+Y>0;\\
&Y+X(nz-(z+1)^2)<0,X-4nz^2<0. \end{split}\] From case (2), we get
the condition (1.5)
\[X<0,-\frac{X(X-2z(nz+(z+1)^2))}{2z}<Y<((z+1)^2-nz)X,\]
where $nz-(z+1)^2>0$. This completes the proof of Theorem 1.1.  \hfill $\oblong$\\

\textbf{Example 1.} When $n=17,$ from $17z-(z+1)^2>0$, we have
\[0.06696562634 \approx\frac{15}{2}-\frac{\sqrt{221}}{2}<z<\frac{15}{2}+\frac{\sqrt{221}}{2}\approx14.93303437.\]
Taking $z=1,$ the curve $E_{17,1}$ has a rational point
$(X,Y)=(-16,-16)$ satisfying the condition (1.5), which leads to the
positive rational solution
$(x,y,z)=\bigg(\frac{4}{7},\frac{2}{3},1\bigg)$. Then for $n=17$,
Eq. (1.3) has a positive integer solution $(x,y,z,w)=(12,14,21,21)$.

When $n=17,z=3,$ we don't find a rational point on $E_{17,3}$
satisfying the condition (1.5).

\begin{figure}[h!]
  \centering
  \subfigure[]{
    \label{fig:subfig:a} 
    \includegraphics[width=2.2in]{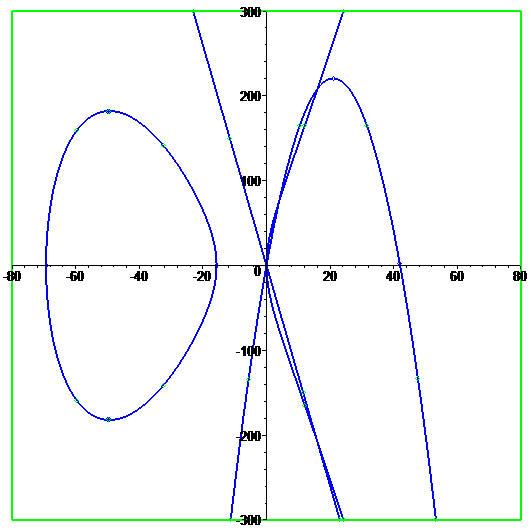}}
  \hspace{0.12in}
  \subfigure[]{
    \label{fig:subfig:b} 
    \includegraphics[width=2.2in]{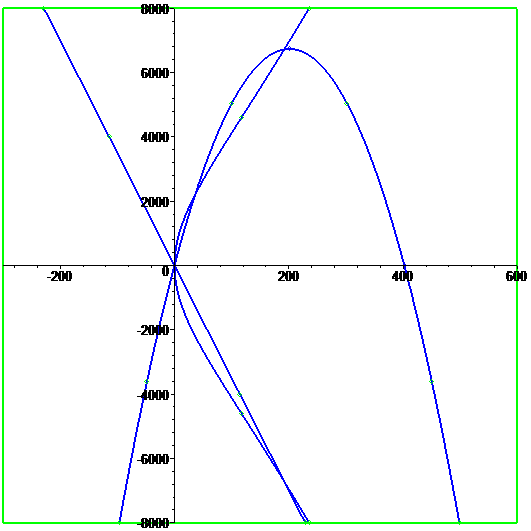}}
  \caption{(a)~$n=17,z=1$ and (b)~$n=17,z=3$ with condition (1.5)}
  \label{fig:subfig} 
\end{figure}

In figure 1, we can see that if the elliptic curve $E_{17,z}$ has
two components, an ``egg" for $X<0$ and an infinite branch for
$X>0,$ then the rational points on $E_{17,z}$ satisfying (1.5) lie
on the egg.

\vskip10pt
\section{Some related questions}
Since Theorem 1.1 is conditional, we pose the following question

\begin{question}
Whether the elliptic curve $E_{n,z}$ has a rational point satisfying
the condition (1.5) for every $n>16$ and some $z>0$ with
$nz-(z+1)^2>0$?
\end{question}

In Table 1 we list a positive integer solution of Eq. (1.3) for
$17\leq n\leq100$ but $n=36,40,64,68,100$ in the range $1\leq
x\leq500, x\leq y \leq 3000, y\leq z\leq 6000$. Fix $w=1$, for some
small $z$, we cannot find a suitable rational point on $E_{n,z}$
such that Eq. (1.4) has a positive rational solution for
$n=36,40,64,68,100$.

\[\begin{tabular}{|c|c||c|c||c|c|}
\hline
 $n$ & $(x,y,z,w)$ & $n$ & $(x,y,z,w)$ &  $n$ & $(x,y,z,w)$ \\
\hline $17$ & $(2, 3, 3, 4)$  &$45$ & $(1, 1, 6, 12)$ &$73$ & $(5,44,45,198)$ \\
\hline $18$ & $(1, 1, 2, 2)$  &$46$ & $(6, 35, 78, 91)$ &$74$ & $(28,33,209,756)$\\
\hline $19$ & $(5, 8, 12, 15)$ &$47$ & $(7, 12, 78, 91)$ &$75$ & $(4, 7, 78, 91)$\\
\hline $20$ & $(1, 1, 1, 3)$  &$48$ & $(1, 1, 3, 15)$ &$76$ & $(1, 7, 10, 42)$\\
\hline $21$ & $(8, 14, 15, 35)$ &$49$ & $(1, 2, 5, 20)$ &$77$ & $(1, 5, 18, 36)$\\
\hline $22$ & $(1, 1, 2, 4)$ &$50$ & $(1, 2, 9, 18)$ &$78$ & $(1, 6, 28, 28)$\\
\hline $23$ & $(76,220,285,385)$ &$51$ & $(35,77,480,528)$ &$79$ & $(1, 3, 24, 28)$\\
\hline $24$ & $(1, 2, 3, 6)$ &$52$ & $(1, 3, 4, 24)$ &$80$ & $(1, 5, 9, 45)$\\
\hline $25$ & $(1, 1, 4, 4)$ &$53$ & $(2, 4, 9, 45)$ &$81$ & $(25,30,65,780)$\\
\hline $26$ & $(20,27,39,130)$ &$54$ & $(1, 3, 8, 24)$ &$82$ & $(7, 24, 112, 273)$\\
\hline $27$ & $(3, 7, 8, 24)$ &$55$ & $(27,132,231,702)$ &$83$ & $(8, 78, 129, 344)$\\
\hline $28$ & $(2, 9, 10, 15)$ &$56$ & $(32,60,207,736)$ &$84$ & $(1, 3, 5, 45)$\\
\hline $29$ & $(1, 1, 4, 6)$ &$57$ & $(3, 6, 40, 56)$ &$85$ & $(1, 18, 20, 36)$\\
\hline $30$ & $(2, 3, 10, 15)$ &$58$ & $(2, 11, 20, 55)$ &$86$ & $(5, 28, 30, 252)$ \\
\hline $31$ & $(1, 4, 5, 10)$ &$59$ & $(25,41,72,600)$ &$87$ & $(2, 4, 15, 84)$ \\
\hline $32$ & $(1, 2, 6, 9)$ &$60$ & $(2, 21, 35, 42)$ &$88$ & $(2, 9, 22, 99)$ \\
\hline $33$ & $(12,35,51,140)$ &$61$ & $(2, 7, 15, 60)$ &$89$ & $(1, 1, 12, 28)$ \\
\hline $34$ & $(6, 35, 40, 63)$ &$62$ & $(3, 16, 45, 80)$ &$90$ & $(3, 21, 80, 120)$ \\
\hline $35$ & $(8, 45, 63, 84)$ &$63$ & $(3, 12, 50, 75)$ &$91$ & $(18, 93, 308, 868)$ \\
\hline $36$ & $(?)$ &$64$ & $(?)$ &$92$ & $(1, 3, 12, 48)$\\
\hline $37$ & $(1, 3, 8, 12)$ &$65$ & $(2, 9, 44, 44)$ &$93$ & $(3, 7, 30, 140)$\\
\hline $38$ & $(2, 3, 15, 20)$ &$66$ & $(1, 12, 12, 30)$ &$94$ & $(1, 5, 8, 56)$\\
\hline $39$ & $(4, 18, 20, 63)$ &$67$ & $(1, 4, 20, 25)$ &$95$ & $(3, 8, 88, 99)$\\
\hline $40$ & $(?)$ &$68$ & $(?)$ &$96$ & $(1, 7, 30, 42)$\\
\hline $41$ & $(1, 5, 12, 12)$ &$69$ & $(24,140,561,595)$ &$97$ & $(5, 20, 21, 276)$\\
\hline $42$ & $(1, 1, 4, 12)$ &$70$ & $(1, 6, 21, 28)$ &$98$ & $(1, 18, 33, 36)$\\
\hline $43$ & $(5, 14, 44, 77)$ &$71$ & $(1, 10, 21, 28)$ &$99$ & $(1, 4, 20, 50)$\\
\hline $44$ & $(2, 14, 15, 35)$ &$72$ & $(1, 4, 21, 28)$ &$100$ & $(?)$\\
\hline
\end{tabular}
\]
\begin{center}Table~~1. A positive integer solution of Eq. (1.3) for $17\leq n\leq100$ \end{center}

\begin{question}
Are there infinitely many positive integers $n$ can be representable
as the product of the sum of four \textbf{distinct positive}
integers with the sum of their reciprocals.
\end{question}

When $m=5,$ by some calculations, for $25\leq n\leq 200$, Eq. (1.1)
has positive integer solutions. So we raise

\begin{question}
For $m \geq 5$, does Eq. (1.1) have a positive integer solution
$(x_1,x_2,\cdot \cdot \cdot,x_m)$ for every integer $n\geq m^2$?
\end{question}

Moreover, we ask
\begin{question}
Whether $m=5$ is the least number such that the integers $n\geq m^2$
can be representable as the product of the sum of $m$ positive
integers with the sum of their reciprocals?
\end{question}

When $m=5,$ we have
\[\begin{split}
&36=(1+1+2+4+4)\bigg(1+1+\frac{1}{2}+\frac{1}{4}+\frac{1}{4}\bigg),\\
&40=(2+9+9+10+15)\bigg(\frac{1}{2}+\frac{1}{9}+\frac{1}{9}+\frac{1}{10}+\frac{1}{15}\bigg),\\
&64=(1+2+3+4+20)\bigg(1+\frac{1}{2}+\frac{1}{3}+\frac{1}{4}+\frac{1}{20}\bigg),\\
&68=(1+5+12+15+15)\bigg(1+\frac{1}{5}+\frac{1}{12}+\frac{1}{15}+\frac{1}{15}\bigg),\\
&100=(1+1+10+15+18)\bigg(1+1+\frac{1}{10}+\frac{1}{15}+\frac{1}{18}\bigg).
\end{split}\]

\begin{remark} In $1963$, R.L. Graham \cite{Graham} proved that every integer greater
than $77$ can be partitioned into \textbf{distinct} positive
integers whose reciprocals add to $1$, i.e.,
\[n=x_1+x_2+\cdot\cdot\cdot+x_m,~\frac{1}{x_1}+\frac{1}{x_2}+\cdot\cdot\cdot+\frac{1}{x_m}=1,\]
where $0<x_1< x_2< \cdot\cdot\cdot< x_m.$ Then for $n>77$ there
exist \textbf{distinct} positive integers
$x_1,x_2,\cdot\cdot\cdot,x_m$ such that
\[n=(x_1+x_2+\cdot\cdot\cdot+x_m)\bigg(\frac{1}{x_1}+\frac{1}{x_2}+\cdot\cdot\cdot+\frac{1}{x_m}\bigg).\]

In $2013$, G. K\"{o}hler and J. Spilker \cite{Kohler-Spilker} showed
that every integer greater than $23$ can be partitioned into
positive integers whose reciprocals add to $1$ (which is called
harmonic partitions). Then the integers $n>23$ can be representable
as the product of the sum of $m$ positive integers with the sum of
their reciprocals for some $m$.

\end{remark}

\vskip20pt
\bibliographystyle{amsplain}

\end{document}